\documentclass[a4paper,12pt]{amsart}
\usepackage{amssymb,amscd,amsmath,a4wide,graphicx,stmaryrd,fullpage,setspace,microtype,xr,textcomp,enumitem,wasysym}
\usepackage[pagebackref=true]{hyperref}
\usepackage[all]{xy}
\usepackage[T1]{fontenc}
\usepackage{lmodern}
\usepackage[numbers,sort&compress]{natbib}

\title{Comparing the non-unital and unital settings for directed homotopy}

\author[P. Gaucher]{Philippe Gaucher}

\address{Universit\'e Paris Cit\'e, CNRS, IRIF, F-75013, Paris, France}

\urladdr{http://www.irif.fr/{\~{}}gaucher} 

\makeatletter
\@namedef{subjclassname@2020}{\textup{2020} Mathematics Subject Classification}
\makeatother


\swapnumbers

\newcommand{\C}{\mathcal{C}}
\newcommand{\D}{\mathcal{D}}
\newcommand{\K}{\mathcal{K}}

\newcommand{\W}{\mathcal{W}}
\newcommand{\F}{\mathcal{F}}

\newcommand{\p}{\times}
\renewcommand{\vec}{\overrightarrow}
\renewcommand{\P}{\mathbb{P}}

\newtheorem*{thmN}{Theorem}

\newtheorem{thm}{Theorem}[section]
\newtheorem{prop}[thm]{Proposition}

\newtheorem{lem}[thm]{Lemma}
\newtheorem{exa}[thm]{Example}
\newtheorem{cor}[thm]{Corollary}

\newtheorem{defn}[thm]{Definition}
\newtheorem{nota}[thm]{Notation}

\newcommand{\bd}{\begin{defn}}
	\newcommand{\ed}{\end{defn}}
\newcommand{\bp}{\begin{prop}}
	\newcommand{\ep}{\end{prop}}
\newcommand{\bth}{\begin{thm}}
	\renewcommand{\eth}{\end{thm}}
\newcommand{\bpf}{\begin{proof}}
	\newcommand{\epf}{\end{proof}}
\newcommand{\bc}{\begin{cor}}
	\newcommand{\ec}{\end{cor}}

\newcommand{\fL}[1]{\ar@{->}[ll]_-{#1}}
\newcommand{\fR}[1]{\ar@{->}[rr]^-{#1}}
\newcommand{\fRr}[1]{\ar@{->}[rrr]^-{#1}}
\newcommand{\fD}[1]{\ar@{->}[dd]_-{#1}}
\newcommand{\fU}[1]{\ar@{->}[uu]^-{#1}}
\newcommand{\f}[2]{\ar@{->}[#1]|{#2}}
\newcommand{\ff}[2]{\ar@2{->}[#1]|{#2}}
\newcommand{\frr}[1]{\ar@{->}[rrrr]^-{#1}}

\newcommand{\fl}[1]{\ar@{->}[l]_-{#1}}
\newcommand{\fr}[1]{\ar@{->}[r]^-{#1}}
\newcommand{\fd}[1]{\ar@{->}[d]_-{#1}}
\newcommand{\fu}[1]{\ar@{->}[u]^-{#1}}

\newcommand{\uu}{{\mathbf{\underline{3}}}}
\renewcommand{\top}{{\mathbf{Top}}}

\newcommand{\iso}{\cong}

\DeclareMathOperator{\Mon}{Mon}

\newcommand{\vI}{\vec{I}}
\renewcommand{\leq}{\leqslant}
\renewcommand{\geq}{\geqslant}

\def\cartesien{%
	\ar@{-}[]+R+<6pt,-2pt>;[]+RD+<6pt,-6pt>%
	\ar@{-}[]+D+<2pt,-6pt>;[]+RD+<6pt,-6pt>%
}
\def\cocartesien{%
	\ar@{-}[]+L+<-6pt,+2pt>;[]+LU+<-6pt,+6pt>%
	\ar@{-}[]+U+<-2pt,+6pt>;[]+LU+<-6pt,+6pt>%
}
\def\hocartesien{%
	\ar@{-}[]+R+<6pt,-2pt>;[]+RD+<6pt,-6pt>_{h}%
	\ar@{-}[]+D+<2pt,-6pt>;[]+RD+<6pt,-6pt>%
}
\def\hococartesien{%
	\ar@{-}[]+L+<-6pt,+2pt>;[]+LU+<-6pt,+6pt>_{h}%
	\ar@{-}[]+U+<-2pt,+6pt>;[]+LU+<-6pt,+6pt>%
}

\newcommand{\brm}[1]{{\rm{\mathbf{#1}}}}

\newcommand{\dtop}{{\brm{Flow}}}

\newcommand{\set}{{\brm{Set}}}

\newcommand{\ttop}{{\brm{TOP}}}

\newcommand{\glob}{{\mathrm{Glob}}}

\newcommand{\ttt}{two-out-of-three property}

\DeclareMathOperator{\id}{Id}

\DeclareMathOperator{\Obj}{Obj}

\newcommand{\mins}[1]{\texorpdfstring{$#1$}{Lg}}

\newcommand{\cat}{{\mathbf{Cat}}}

\newcommand{\cattop}{{\mathbf{Cat}}}

\renewcommand{\P}{\mathbb{P}}

\makeatletter
\def\varholim@#1#2{%
	\vtop{\m@th\ialign{##\cr
			\hfil$#1\operator@font holim$\hfil\cr
			\noalign{\nointerlineskip\kern1.5\ex@}#2\cr
			\noalign{\nointerlineskip\kern-\ex@}\cr}}%
}
\def\holimproj{%
	\mathop{\mathpalette\varholim@{\leftarrowfill@\textstyle}}\nmlimits@
}
\def\holiminj{%
	\mathop{\mathpalette\varholim@{\rightarrowfill@\textstyle}}\nmlimits@
}
\makeatother

\DeclareMathOperator{\cell}{{\brm{cell}}}
\DeclareMathOperator{\cof}{{\brm{cof}}}
\DeclareMathOperator{\inj}{{\brm{inj}}}

\DeclareMathOperator{\cocyl}{{Path}}
\newcommand{\adj}[5]{\xymatrix@C=#5em{{#1}\ar@/^0.8em/[r]^-{#2} \ar@{}[r]|-{\perp} & \ar@/^0.8em/[l]^-{#3} {#4}}}

\DeclareMathOperator{\II}{I^+}

\newcommand{\Di}{\mathbf{D}}
\newcommand{\Sp}{\mathbf{S}}

\begin{document}

\begin{abstract}
	This note explores the link between the q-model structure of flows and the Ilias model structure of topologically enriched small categories. Both have weak equivalences which induce equivalences of fundamental (semi)categories. The Ilias model structure cannot be left-lifted along the left adjoint adding identity maps. The minimal model structure on flows having as cofibrations the left-lifting of the cofibrations of the Ilias model structure has a homotopy category equal to the $3$-element totally ordered set. The q-model structure of flows can be right-lifted to a q-model structure of topologically enriched small categories which is minimal and such that the weak equivalences induce equivalences of fundamental categories. The identity functor of topologically enriched small categories is neither a left Quillen adjoint nor a right Quillen adjoint between the q-model structure and the Ilias model structure.
\end{abstract}
	
\maketitle
\tableofcontents

\section{Introduction}

\subsection*{Presentation}

The time flow of a concurrent process can be modelled by a topologically enriched small semicategory \cite{model3} or by a topologically enriched small category \cite{UnitalDihomotopy,DHH}. The objects represent the \textit{states} of the concurrent process and the nonidentity morphisms represent the \textit{execution paths}, the topology modelling concurrency \cite{DAT_book}. The primary reason for excluding identity morphisms in \cite[Definition~4.11]{model3} is to obtain \textit{functorial} constructions for the branching and merging homology theories (see \cite[Section~20]{model3}). It enables us to prove the invariance by refinement of observation in \cite[Corollary~11.3]{3eme}, and therefore to fix  Goubault-Jensen's construction of \cite{Prehistoire-homology}. The main technical tool is the minimal model category introduced in \cite{model3}, called the \textit{q-model structure (of flows)} after \cite[Theorem~7.6]{QHMmodel}. The examples coming from computer science are non-unital as well because they are modelled by \textit{precubical} sets (e.g. \cite{zbMATH06546404,symcub,ccsprecub,zbMATH07226006}) and because precubical sets have non-unital geometric realizations \cite[Definition~7.2]{ccsprecub}. The transverse degeneracy maps of precubical sets, introduced for the \textit{functorial} formalization of the parallel product with synchronization of process algebra \cite[Theorem~3.1.15 and Definition~4.2.2]{symcub}, belong to the non-unital world as well. The transverse degeneracy maps lead to a vast generalization of Raussen's notion of natural $d$-path in \cite{directed_degeneracy}. The non-unital setting is also necessary to construct the underlying homotopy type functor which is geometrically the homotopy type of the space obtained after forgetting the temporal information  \cite[Section~6]{4eme} \cite[Proposition~8.16]{Moore2}. 

On the other hand, the mathematical literature provides several constructions of model structures on enriched small categories such that the weak equivalences are the so-called Dwyer-Kan equivalences of \cite{Dwyer1980}: for simplicially enriched small categories \cite{Bergner2006}, for topologically enriched small categories \cite{cattopBK} and for small categories enriched in a given monoidal model category \cite{MR3094501}. The generating cofibrations of the q-model structure of flows of \cite{model3} are almost those obtained by transfer along the left adjoint formally adding identity maps from the generating cofibrations of the Ilias model structure constructed in \cite{cattopBK}. The only difference is the presence of the flow cofibration $R:\{0,1\} \to \{0\}$ which has no counterpart in the Ilias model structure (see Proposition~\ref{cof-left-induite}). This leads to the question of comparing the model structures on flows and on topologically enriched small categories. The following sequence of theorems answers the question. 

\begin{thmN} (Theorem~\ref{noleftinduced})
	The Ilias model structure on topologically enriched small categories \cite{cattopBK} cannot be transferred to the category of flows along the left adjoint formally adding identity maps.
\end{thmN}

\begin{thmN} (Corollary~\ref{DK-explicit1})
	The minimal model structure on flows with respect to the transfer of the cofibrations of the Ilias model structure along the left adjoint formally adding identity maps has three homotopy types. 
\end{thmN}

\begin{thmN} (Theorem~\ref{three2})
	The q-model structure of flows can be transferred along the right adjoint forgetting the identity maps to the category of topologically enriched small categories. We obtain a combinatorial model structure which is minimal. Its weak equivalences induce equivalences of fundamental categories. The left Quillen adjoint formally adding identity maps from flows to enriched small categories is not a left Quillen equivalence. 
\end{thmN}

The model category of Theorem~\ref{three2} on topologically enriched small categories seems to be new. With the same argument, the h-model structure and the m-model structure of flows constructed in \cite[Theorem~7.4]{QHMmodel} can be transferred along the right adjoint forgetting the identity maps to the category of topologically enriched small categories as well. We obtain a h-model structure and a m-model structure on topologically enriched small categories which are both accessible as model categories. 

The following table of minimal model categories summarizes the results of this note. The symbol $\CIRCLE$ means that the weak equivalences induce equivalences of fundamental (semi)categories. The symbol $\Circle$ means that they do not. 

\begin{center}
	\begin{tabular}{|c||cc|cc|}
		\hline
		& &$R$ is a cofibration &  &$R$ is not a cofibration \\
		\hline
		\hline
		$\dtop$ & $\CIRCLE$ &q-model structure of \cite{model3} & $\Circle$ &Corollary~\ref{DK-explicit1}\\
		\hline
		$\cattop$ & $\CIRCLE$ &Theorem~\ref{three2} & $\CIRCLE$ & Ilias model structure of \cite{cattopBK} \\
		\hline
	\end{tabular}
\end{center}

The conclusion that must be drawn from this note is that the flow cofibration $R:\{0,1\} \to \{0\}$ is much more important in a globular approach of directed homotopy than what was expected in \cite{model3}.

\subsection*{Prerequisites and notations}

We refer to \cite{TheBook} for locally presentable categories, to \cite{MR2506258,HKRS17,GKR18} for combinatorial and accessible model categories. We refer to \cite{MR99h:55031,ref_model2} for more general model categories. We work with a locally presentable convenient category of topological spaces $\top$ for doing algebraic topology. The internal hom is denoted by $\ttop(-,-)$. The category of \textit{$\Delta$-generated spaces} or of \textit{$\Delta$-Hausdorff $\Delta$-generated spaces} (cf. \cite[Section~2 and Appendix~B]{leftproperflow}) are two such examples. The category $\top$ is equipped with its q-model structure (we use the terminology of \cite{ParamHomTtheory}). What follows is some notations and conventions: $\varnothing$ is the initial object, $\mathbf{1}$ is the final object, $\id_X$ is the identity of $X$. A model structure $(\C,\F,\W)$ means that the class of cofibrations is $\C$, that the class of fibrations is $\F$ and that the class of weak equivalences is $\W$. A combinatorial model structure on $\K$ is \textit{minimal}  if the class of weak equivalences is the smallest Grothendieck localizer with respect to its set of generating cofibrations \cite{rotho,henry2020minimal}. Note that in \cite{rotho}, the adjective \textit{left-determined} is used instead. When all objects of a model category are fibrant, any Grothendieck localizer which is strictly smaller than the class of weak equivalences never induces a model structure. By \cite[Theorem~1.4]{henry2020minimal}, every tractable combinatorial model category with fibrant objects is minimal. The notation $f\boxslash g$ means that $g$ satisfies the right lifting property (RLP) with respect to $f$; ${^\boxslash}\C = \{g, \forall f \in \C, g\boxslash f\}$; $\C^\boxslash = \inj(\C) = \{g, \forall f \in \C, f\boxslash g\}$; $\cof(\C)={^\boxslash}(\C^\boxslash)$; $\cell(\C)$ is the class of transfinite compositions of pushouts of elements of $\C$. A \textit{cellular} object $X$ of a combinatorial model category is an object such that the canonical map $\varnothing\to X$ belongs to $\cell(I)$ where $I$ is the set of generating cofibrations.

In this paper, the transfer of a model structure of $\dtop$ along the right adjoint $\cattop\subset \dtop$ of Proposition~\ref{well-known} is called the \textit{right-lifting} of the model structure of $\dtop$. Similarly, the transfer of a model structure of $\cattop$ (of a weak factorization system resp.) along the left adjoint $\II:\dtop\to \cattop$ of Proposition~\ref{well-known} is called the \textit{left-lifting} of the model structure of $\cattop$ (of the weak factorization system resp.). See the introductions of \cite{HKRS17,GKR18} and the beginning of \cite[Section~2]{GKR18} for further explanations.

\subsection*{Acknowledgments}

I thank Simon Henry for a useful discussion about semicategories. I thank the anonymous referee for the suggestions to improve and clarify the presentation.

\section{The adjunction \mins{\II:\dtop\leftrightarrows \cat:\supset}}

\bd \cite[Definition~4.11]{model3} \label{def-flow}
A {\rm flow} is a small semicategory enriched over the closed monoidal category $(\top,\p)$. The corresponding category is denoted by $\dtop$. 
\ed

A \textit{flow} $X$ consists of a topological space $\P X$ of \textit{execution paths}, a discrete space $X^0$ of \textit{states}, two continuous maps $s$ and $t$ from $\P X$ to $X^0$ called the source and target map respectively, and a continuous and associative map $*:\{(x,y)\in \P X\p \P X; t(x)=s(y)\}\longrightarrow \P X$ such that $s(x*y)=s(x)$ and $t(x*y)=t(y)$. Let $\P_{\alpha,\beta}X = \{x\in \P X\mid s(x)=\alpha \hbox{ and } t(x)=\beta\}$. Note that the composition is denoted by $x*y$, not by $y\circ x$.

Every set can be viewed as a flow with an empty space of execution paths. The obvious functor $\set \subset \dtop$ from the category of sets to that of flows is limit-preserving and colimit-preserving. The following examples of flows are important for the sequel: 

\begin{exa}
	For a topological space $Z$, let $\glob(Z)$ be the flow defined by 
	\[
	\glob(Z)^0=\{0,1\}, \ 
	\P \glob(Z)= \P_{0,1} \glob(Z)=Z,\ 
	s=0,\  t=1.
	\]
	This flow has no composition law. The {\rm directed segment} is the flow $\vI = \glob(\{0\})$.
\end{exa}

\begin{exa} \label{B}
	Denote by $B$ (like branching) the flow $1\leftarrow 0 \rightarrow 1$ with three states and two execution paths. This flow has no composition law.
\end{exa}

\begin{nota}
	Let $n\geq 1$. Denote by $\mathbf{D}^n = \{b\in \mathbb{R}^n, |b| \leq 1\}$ the $n$-dimensional disk, and by $\mathbf{S}^{n-1} = \{b\in \mathbb{R}^n, |b| = 1\}$ the $(n-1)$-dimensional sphere. By convention, let $\mathbf{D}^{0}=\{0\}$ and $\mathbf{S}^{-1}=\varnothing$.
\end{nota}

\begin{nota} Let
	\begin{align*}
	&I^{gl} = \{c_n:\glob(\mathbf{S}^{n-1})\subset \glob(\mathbf{D}^{n}) \mid n\geq 0\},\\
	&J^{gl} = \{\glob(\mathbf{D}^{n}\p\{0\})\subset \glob(\mathbf{D}^{n}\p [0,1] \mid n\geq 0\},\\
	&C:\varnothing \to \{0\},R:\{0,1\} \to \{0\}.
	\end{align*}
\end{nota}

\begin{nota}
	The category of small categories enriched over $\top$ is denoted by $\cattop$. The set of objects of an enriched small category $X$ is denoted by $\Obj(X)$ and the space of morphisms from $A$ to $B$ by $X(A,B)$. 
\end{nota}

\bp (well-known) \label{well-known}
The inclusion $\cattop \subset \dtop$ has a left adjoint~\footnote{It has also a right adjoint, the enriched small category of idempotents of a flow, which is not used in this note.} denoted by $\II:\dtop\longrightarrow \cattop$. It consists of adding identity maps as isolated points in the spaces of morphisms. This functor is faithful.
\ep

What follows is an adaptation of \cite[Definition~4.37]{DAT_book}.

\bd
Let $X$ be an object of $\dtop$. The {\rm fundamental semicategory} of $X$ is the small semicategory $\vec{\pi}_1(X)$ having $X^0$ for the set of objects and the set of morphisms between two objects is the set of path-connected components of the space of execution paths between these two objects. If $X\in \cat\subset \dtop$, then $\vec{\pi}_1(X)$ is a small category which is called the {\rm fundamental category} of $X$.
\ed

For all $X\in \dtop$, $\II(\vec{\pi}_1(X))$ is also a small category which is called the fundamental category of $X$. For $X\in \cat$, the canonical map $\II(\vec{\pi}_1(X)) \to \vec{\pi}_1(X)$ is not an equivalence of categories. 

\section{Left-lifting the Ilias model structure}

\bth \cite[Theorem~2.4]{cattopBK} \label{DK}
There exists one and only one combinatorial model structure $(\cat)_{DK} = (\C_{DK},\F_{DK},\W_{DK})$ on $\cattop$ with the following properties: 
\begin{itemize}[leftmargin=*]
	\item A set of generating cofibrations is the set of maps $\II(I^{gl}\cup\{C\})$.
	\item The weak equivalences are the DK-equivalences: there are the maps of enriched functors $F:\C\to \D$ such that ${\vec{\pi}_1}(F):{\vec{\pi}_1}(\C) \to {\vec{\pi}_1}(\D)$ is an equivalence of categories and such that for all pairs of objects $(\alpha,\beta)$ of $\C$, there is a weak homotopy equivalence $\C(\alpha,\beta)\to \D(F(\alpha),F(\beta))$.
	\item A set of generating trivial cofibrations is given by the set of maps $\II(J^{gl})\cup \{\II(\{0\}) \to (\{0\iso 1\})^{cof}\}$ where $\{0\iso 1\}$ is the small category with two isomorphic objects $0$ and $1$.
\end{itemize}
It is called the {\rm Ilias model structure}. All objects are fibrant.
\eth

Theorem~\ref{DK} is the topological analogue of the Bergner model structure on simplicially enriched small categories \cite{Bergner2006}. The weak equivalences are the \textit{Dwyer-Kan equivalences} of \cite{Dwyer1980}. The combinatorial model category is minimal since all objects are fibrant. The weak equivalences of $\W_{DK}$ induce equivalences of fundamental categories by definition.

\bp \label{inv-square}
Let $f:X\to Y$ be a map of flows. Let $i:A\to B\in I^{gl}\cup\{C\}$. Consider a commutative square of $\cattop$
\[
\xymatrix@C=3em@R=3em
{
	\II(A) \fr{\phi} \fd{\II(i)} & \II(X) \ar@{->}[d]^-{\II(f)} \\
	\II(B) \fr{\overline{\phi}} &  \II(Y)
}
\]
Then either $\II(B)\sqcup_{\II(A)}\II(X) \iso \II(X)$, or the canonical map $\II(B)\sqcup_{\II(A)}\II(X) \to \II(Y)$ is of the form $\II(g)$ for some unique map of flows $g:B\sqcup_{A}X \to Y$. 
\ep
 

\bpf 
That there is at most one such a map $g$ is a consequence of the fact that $\II$ is faithful. A commutative diagram of enriched small categories of the form 
\[
\xymatrix@C=3em@R=3em
{
	\varnothing \fr{} \fd{} & \II(X) \ar@{->}[d]^-{\II(f)} \\
	\II(\{0\}) \fr{} &  \II(Y)
}
\]
is the image by the functor $\II:\dtop\to \cattop$ of the commutative diagram of flows 
\[
\xymatrix@C=3em@R=3em
{
	\varnothing \fr{} \fd{} & X \ar@{->}[d]^-{f} \\
	\{0\} \fr{} &  Y
}
\]
Thus, in this case, $g$ exists by the universal property of the pushout. Consider now a commutative diagram $(C)$ of enriched small categories of the form 
\[
\xymatrix@C=3em@R=3em
{
	\II(\glob(\mathbf{S}^{n-1})) \fr{\phi} \fd{} & \II(X) \ar@{->}[d]^-{\II(f)} \\
	\II(\glob(\mathbf{D}^{n})) \fr{\overline{\phi}} &  \II(Y)
}
\]
with $n\geq 0$. If $\phi(0)\neq \phi(1)$ are two different objects of $\II(X)$, then the latter commutative diagram of enriched small categories is the image by the functor $\II:\dtop\to \cattop$ of the commutative diagram $(D)$ of flows
\[
\xymatrix@C=3em@R=3em
{
	\glob(\mathbf{S}^{n-1}) \fr{\phi} \fd{} & X \ar@{->}[d]^-{f} \\
	\glob(\mathbf{D}^{n}) \fr{} &  Y
}
\]
We conclude the existence of $g$ as above. It remains the case $\phi(0)=\phi(1)$. In this case, we have the commutative diagram of topological spaces
\[
\xymatrix@C=3em@R=3em
{
	\mathbf{S}^{n-1} \fr{\phi} \fd{} & \{\id_{\phi(0)}\} \sqcup \P_{\phi(0),\phi(1)}X \ar@{->}[d]^-{f} \\
	\mathbf{D}^{n} \fr{\overline{\phi}} &  \{\id_{\phi(0)}\} \sqcup \P_{\phi(0),\phi(1)}Y
}
\]
If $n\geq 1$ and since $\mathbf{D}^{n}$ is connected, then either $\overline{\phi}(\mathbf{D}^{n}) \subset \{\id_{\phi(0)}\}$ and $\phi(\mathbf{S}^{n-1})\subset \{\id_{\phi(0)}\}$ or $\overline{\phi}(\mathbf{D}^{n}) \subset \P_{\phi(0),\phi(1)}Y$ and $\phi(\mathbf{S}^{n-1})\subset \P_{\phi(0),\phi(1)}X$. If $n=0$, then $\mathbf{S}^{n-1}=\varnothing$ and either $\overline{\phi}(\mathbf{D}^{n}) \subset \{\id_{\phi(0)}\}$ or $\overline{\phi}(\mathbf{D}^{n}) \subset \P_{\phi(0),\phi(1)}Y$. 

In the first alternative in both cases, there is the pushout diagram of enriched small categories 
\[
\xymatrix@C=3em@R=3em
{
	\II(\glob(\mathbf{S}^{n-1})) \fr{\phi} \fd{} & \II(X) \ar@{->}[d]^-{\II(f)} \\
	\II(\glob(\mathbf{D}^{n})) \fr{} &  \cocartesien \II(X).
}
\]

In the second alternative in both cases, the commutative diagram $(C)$ is the image by the functor $\II:\dtop\to \cattop$ of the commutative diagram $(D)$ and we conclude the existence of $g$ as above.
\epf

By \cite[Theorem~2.6]{GKR18}, the left-lifting of the small weak factorization system $(\C_{DK},\F_{DK}\cap \W_{DK})$ along the left adjoint $\II:\dtop\to \cat$ exists and is accessible. In fact, we have the proposition:

\bp \label{cof-left-induite}
The left-lifting of the small weak factorization system $(\C_{DK},\F_{DK}\cap \W_{DK})$ along the left adjoint $\II:\dtop\to \cat$ is small, being generated by $I^{gl}\cup \{C\}$. 
\ep

\bpf
It suffices to prove that $\II^{-1}(\C_{DK}) =\cof(I^{gl}\cup \{C\})$. We have $\II(I^{gl}\cup \{C\}) \subset \C_{DK}$ by Theorem~\ref{DK}. Since $\II:\dtop\to \cattop$ is a left adjoint, we obtain the inclusion $\cell(I^{gl}\cup \{C\}) \subset \II^{-1}(\C_{DK})$. And using the fact that every map of $\cof(I^{gl}\cup \{C\})$ is a retract of a map of $\cell(I^{gl}\cup \{C\})$, we obtain the inclusion $\cof(I^{gl}\cup \{C\}) \subset \II^{-1}(\C_{DK})$ since the class of maps $\C_{DK}$ is closed under retract. Conversely, let $f:X\to Y$ be a map of flows such that $\II(f):\II(X)\to \II(Y)$ is a cofibration of $\cattop$. By using the small object argument of \cite[Theorem~2.1.14]{MR99h:55031}, we factor $\II(f)$ as a composite $\II(X) \to Z \to \II(T)$ such that the map $\II(X) \to Z$ belongs to $\cell(\II(I^{gl}\cup\{C\}))$ and such that the map $Z\to \II(T)$ belongs to $\inj(\II(I^{gl}\cup\{C\}))$. Since $\II$ is a left adjoint, by an immediate transfinite induction, there exists a transfinite tower $(X_\alpha)_{\alpha<\lambda}$ of $\dtop$ with $X=X_0$ and $Z=\II(X_\lambda)$ such that each map $X_\alpha \to X_{\alpha+1}$ for $\alpha<\lambda$ is a pushout of a map of $I^{gl}\cup \{C\}$. By induction on $\alpha\geq 0$, let us prove that the map of enriched small categories $\II(X_\alpha)\to \II(T)$ is the image by the functor $\II$ of a map of flows $g_\alpha:X_\alpha\to T$. There is nothing to prove for $\alpha=0$. The passage from $\alpha$ to $\alpha+1$ is ensured by Proposition~\ref{inv-square}. Finally, the statement holds for a limit ordinal $\alpha$ since $\II$ is colimit-preserving. We deduce that the map of enriched small categories $Z\to \II(T)$ is of the form $\II(g)$ for some map of flows $g:X_\lambda \to T$: take $g=g_\lambda$. The lift $\ell$ in the commutative diagram of enriched small categories
\[
\xymatrix@C=3em@R=3em
{
	\II(X) \fr{} \fd{\II(f)} & \II(X_\lambda) \ar@{->}[d]^-{\II(g)} \\
	\II(T) \ar@{=}[r] \ar@{->}[ru]^-{\ell} &  \II(T)
}
\]
exists since $\II(f)$ is a cofibration of $\cattop$ by hypothesis. For all $\alpha\in T^0$, the commutativity of the diagram of spaces
\[
\xymatrix@C=3em@R=3em
{
	\{\id_\alpha\} \sqcup \P_{\alpha,\alpha} T\ar@{=}@/^25pt/[rr] \fr{\ell} & \{\id_{\ell(\alpha)}\} \sqcup \P_{\ell(\alpha),\ell(\alpha)} X_\lambda \fr{\II(g)} & \{\id_\alpha\} \sqcup \P_{\alpha,\alpha} T
}
\]
implies that $\ell(\P_{\alpha,\alpha} T)\subset \P_{\ell(\alpha),\ell(\alpha)} X_\lambda$, and therefore that $\ell=\II(\overline{\ell})$ for some map of flows $\overline{\ell}:T\to X_\lambda$. Since the functor $\II$ is faithful, we obtain the commutative diagram of flows 
\[
\xymatrix@C=3em@R=3em
{
	X \fr{} \fd{f} & X_\lambda \ar@{->}[d]^-{g} \\
	T \ar@{=}[r] \ar@{->}[ru]^-{\overline{\ell}} &  T
}
\]
It means that the map of flows $f:X\to T$ is a retract of the map of flows $X\to X_\lambda$. We deduce $f\in \cof(I^{gl}\cup \{C\})$, the map $X\to X_\lambda$ belonging to $\cell(I^{gl}\cup \{C\})$ by construction. We deduce the inclusion $\II^{-1}(\C_{DK}) \subset \cof(I^{gl}\cup \{C\})$.
\epf

\begin{lem} \label{isset}
	Let $f:X\to Y$ be a map of $\dtop$ such that $Y$ is a set. Then $X$ is a set as well. 
\end{lem}

\bpf It is a consequence of the lack of identity maps for the objects of $\dtop$.
\epf

\bth \label{noleftinduced}
The model category $(\cattop)_{DK}$ cannot be left-lifted along the left adjoint $\II:\dtop\to \cattop$. 
\eth

\bpf
By Proposition~\ref{cof-left-induite} and Lemma~\ref{isset}, the map $R:\{0,1\}\to \{0\}$ satisfies the RLP with respect to $\II^{-1}(\C_{DK})$ because it satisfies the RLP with respect to $C:\varnothing\to\{0\}$. But $\II(R)\notin \W_{DK}$. It means that the left acyclicity condition $\II^{-1}(\C_{DK})^\boxslash \subset \II^{-1}(\W_{DK})$ fails and that the left-induced model structure does not exist by \cite[Proposition~2.3]{GKR18}. 
\epf

Theorem~\ref{noleftinduced} can be proved without using Proposition~\ref{cof-left-induite}. Indeed, thanks to Lemma~\ref{isset}, the only maps of flows $f$ belonging to $\II^{-1}(\C_{DK})$ such that there exists a morphism in the category of maps of flows from $f$ to $R$ are the set maps of $\cell(C)=\cof(C)$, i.e. the one-to-one set maps. Proposition~\ref{cof-left-induite} is proved because it is used in Corollary~\ref{DK-explicit1}.

\section{Left-lifting the cofibrations of the Ilias model structure}

We need to recall:

\bth \cite[Theorem~7.6]{QHMmodel} \label{q}
There exists one and only one combinatorial model structure $(\dtop)_{q}$ on $\dtop$ with the following properties: 
\begin{itemize}[leftmargin=*]
	\item A set of generating cofibrations is the set of maps $I^{gl}\cup\{C,R\}$.
	\item The weak equivalences are the maps of flows $f:X\to Y$  inducing a bijection $f^0:X^0\iso Y^0$ and a weak homotopy equivalence $\P f:\P X \to \P Y$.
	\item A set of generating trivial cofibrations is given by the set of maps $J^{gl}$.
\end{itemize}
It is called the {\rm q-model structure}. The cofibrations (fibrations resp.) are called q-cofibrations (q-fibrations resp.). All flows are q-fibrant.
\eth

The weak equivalences of $(\dtop)_q$ induce isomorphisms of fundamental semicategories. The q-model structure of flows is minimal by \cite[Theorem~1.4]{henry2020minimal} since it is combinatorial and all its objects are fibrant~\footnote{\cite[Theorem~4.3]{leftdetflow} gives another argument which does not require to use a locally presentable setting.}.

\bd 
The class of maps of flows $\overline{\W}_{DK}$ consists of the maps of flows $f:X\to Y$ such that either $X=Y=\varnothing$, or $X$ and $Y$ are both nonempty sets, or $X$ and $Y$ both contain at least one execution path. 
\ed

As an immediate consequence of the definition above, we obtain:

\bp All maps of flows 
\begin{align*}
&I^{gl}_{\geq 1}=\{c_n\mid n\geq 1\},C^+:\{0\}\subset \{0,1\},\\
&c_0^+=\id_{\vI}\sqcup c_0:\vI\sqcup\glob(\mathbf{S}^{-1})\subset \vI\sqcup\glob(\mathbf{D}^{0})
\end{align*}
belong to $\overline{\W}_{DK}$.
\ep

We recall the four following propositions for the convenience of the reader.

\bp \label{globRLP} \cite[Proposition~13.2]{model3}
Let $i:U\to V$ be a map of $\top$. A morphism of flows $f:X\to Y$ satisfies the RLP with respect to $\glob(i):\glob(U)\to \glob(V)$ if and only if for all $(\alpha,\beta)\in X^0\p X^0$, the map $\P_{\alpha,\beta}X\to \P_{f(\alpha),f(\beta)}Y$ satisfies the RLP with respect to $i$. 
\ep

\bp (\cite[Theorem~2.1.19]{MR99h:55031}) \label{rec-mod}
Let $I$ and $J$ be two sets of maps of a locally presentable category $\K$. Let $\W$ be a class of maps satisfying the \ttt\ and which is closed under retract. If $\cell(J)\subset \W\cap \cof(I)$, $\inj(I)\subset \W\cap \inj(J)$ and $\W\cap \cof(I)\subset \cof(J)$, then $(\cof(I),\inj(J),\W)$ is a model structure on $\K$.
\ep

\bp \label{pseudo-leftproper} (\cite[Lemma~5.2.6]{MR99h:55031})
Let $\mathcal{M}$ be a model category. Consider a pushout diagram of $\mathcal{M}$ of the form
\[
\xymatrix@C=3em@R=3em
{
	X  \ar@{^(->}[r] \fd{\simeq} & Y \fd{} \\
	Z \fr{} & \cocartesien T
}
\]
such that $X,Y,Z$ are cofibrant, such that the top horizontal map is a cofibration and such that the left vertical map is a weak equivalence. Then the right vertical map $Y\to T$ is a weak equivalence.
\ep

\bp \label{connected-colim-glob} \cite[Proposition~3.7]{leftdetflow}
The globe functor $\glob:\top\to \dtop$ preserves connected colimits (i.e. colimits such that the underlying small category is connected).
\ep

\begin{nota}
	Let $\uu$ be the small category associated with the poset $\{0\leq 1\leq 2\}$. 
\end{nota}

\bth \label{DK-explicit}
There exists one and only one model structure on $\dtop$ such that 
\begin{itemize}[leftmargin=*]
	\item A set of generating cofibrations is $I^{gl}\cup \{C\}$.
	\item A set of generating trivial cofibrations is $\{C^+,c_0^+\} \cup J^{gl} \cup I^{gl}_{\geq 1}$.
	\item The class of weak equivalences is $\overline{\W}_{DK}$.
	\item The homotopy category of this model structure is the category $\uu$: every flow is weakly equivalent either to the initial or terminal flow, or to a singleton.
	\item The cofibrant flows are the q-cofibrant flows.
	\item The fibrant flows are the flows $X$ such that $\P X=\varnothing$ (i.e. the sets) and the flows $X$ such that for all $(\alpha,\beta)\in X^0\p X^0$, the space  $\P_{\alpha,\beta}X$ is contractible. In particular, not all flows are fibrant. 
\end{itemize}
Moreover, this combinatorial model structure is minimal.
\eth

\bpf The uniqueness comes from the fact that a model structure is characterized by its cofibrations and its trivial cofibrations. Note that $\uu$ is the full subcategory of $\dtop$ generated by the initial and terminal flows and by the singleton. Consider the functor $\underline{w}:\dtop \to \uu$ characterized as the unique functor which takes a flow $X$ to $0$ if $X^0 = \varnothing$, to $1$ if $X^0\neq \varnothing$ and $\P X=\varnothing$, and to $2$ otherwise. Then $\overline{\W}_{DK}$ is the inverse image by $\underline{w}$ of the identity maps of $\uu$. We deduce that the class $\overline{\W}_{DK}$ has the \ttt\  and that it is closed under retract. 

All maps of $\cell(\{C^+,c_0^+\} \cup J^{gl} \cup I^{gl}_{\geq 1})$ are q-cofibrations which are one-to-one on states, which implies $\cell(\{C^+,c_0^+\} \cup J^{gl} \cup I^{gl}_{\geq 1}) \subset \cof(I^{gl}\cup \{C\})$. Every map of $\cell(\{C^+,c_0^+\} \cup J^{gl} \cup I^{gl}_{\geq 1})$ is either between nonempty sets or between flows containing execution paths, hence the inclusion $\cell(\{C^+,c_0^+\} \cup J^{gl} \cup I^{gl}_{\geq 1}) \subset \overline{\W}_{DK}$. 

We obtain the inclusion $\cell(\{C^+,c_0^+\} \cup J^{gl} \cup I^{gl}_{\geq 1})\subset \overline{\W}_{DK}\cap\cof(I^{gl}\cup \{C\})$.  An element $f:X\to Y$ of $\inj(I^{gl}\cup \{C\})$ is surjective on states. Therefore $X^0=\varnothing$ if and only if $Y^0=\varnothing$ and $f\in \inj(C^+)$. By Proposition~\ref{globRLP}, every map $\P_{\alpha,\beta}X\to \P_{f(\alpha),f(\beta)}Y$ for all $(\alpha,\beta)\in X^0\p X^0$ is a trivial q-fibration of spaces. Consequently, $X$ contains execution paths if and only if $Y$ contains execution paths. We deduce that $f\in \overline{\W}_{DK}$. By Proposition~\ref{globRLP} again, we deduce that $f\in \inj(J^{gl} \cup I^{gl})$. We obtain the inclusions $\inj(I^{gl}\cup \{C\})\subset \overline{\W}_{DK}\cap \inj(\{C^+\} \cup J^{gl} \cup I^{gl})\subset \overline{\W}_{DK}\cap \inj(\{C^+,c_0^+\} \cup J^{gl} \cup I^{gl}_{\geq 1})$. 

Finally, a map $f\in \overline{\W}_{DK}\cap \cof(I^{gl}\cup \{C\})$ is a q-cofibration which is one-to-one on states such that either the source and the target are empty, or the source and the target are nonempty set (in this case, $f$ belongs to $\cof(\{C^+\})$), or such that both the source and the target contain execution paths. In the latter case, it belongs to $\cof(\{c_0^+\} \cup I^{gl}_{\geq 1})$. We deduce that $\overline{\W}_{DK}\cap \cof(I^{gl}\cup \{C\})\subset \cof(\{C^+,c_0^+\} \cup J^{gl} \cup I^{gl}_{\geq 1})$. 

The proof of the existence of the model structure is complete thanks to Proposition~\ref{rec-mod}. 

Since all flows are q-fibrant, a flow $X$ is fibrant if and only if the canonical map $X\to \mathbf{1}$ satisfies the RLP with respect to $\{C^+,c_0^+\} \cup I^{gl}_{\geq 1}$. Since $\inj(C^+) \cap \set$ is equal to the surjective set maps union the set maps starting from the empty set by \cite[Lemme~4.4(3)]{nonexistence}, the canonical map $X\to \mathbf{1}$ always satisfies the RLP with respect to $C^+$. Thus a flow $X$ is fibrant if and only if the canonical map $X\to \mathbf{1}$ satisfies the RLP with respect to $\{c_0^+\} \cup I^{gl}_{\geq 1}$. We deduce that all sets viewed as flows are fibrant. Consider now a flow $X$ such that $\P X\neq \varnothing$. Then the map $X\to \mathbf{1}$ satisfies the RLP with respect to $c_0^+$ if and only it satisfies the RLP with respect to $c_0$. The characterization of fibrant objects is complete thanks to Proposition~\ref{globRLP}.

Since not all flows are fibrant for this model structure, an additional argument is required to prove that it is indeed minimal. 

Consider a model structure $(\C,\F,\W)$ on $\dtop$ such that $\C=\cof(I^{gl}\cup \{C\})$. The cofibrant flows are the q-cofibrant flows and the cofibrations are the q-cofibrations which are one-to-one on states. All trivial q-fibrations are trivial fibrations since they satisfy the RLP with respect to $I^{gl}\cup \{C\} \subset I^{gl}\cup \{C,R\}$. 

Observe at first that $R:\{0,1\} \to \{0\}$ is a trivial fibration. We have $R.C^+ = \id_{\{0\}}$. By the \ttt, we deduce that $C^+:\{0\}\subset \{0,1\}$ is a weak equivalence. It means that two nonempty sets viewed as flows are always weakly equivalent. 

We are going to prove by induction on $n\geq 1$ that the map $c_n:\glob(\mathbf{S}^{n-1})\subset \glob(\mathbf{D}^n)$ is a trivial cofibration. From the pushout diagram (see Example~\ref{B})
\[
\xymatrix@C=3em@R=3em
{
	\{1\}\sqcup \{1\} \fr{} \fd{} & B \fd{} \\
	\{1\} \fr{} & \cocartesien \glob(\mathbf{S}^0)
}
\]
and Proposition~\ref{pseudo-leftproper}, we deduce that the map $B\to \glob(\mathbf{S}^0)$ is a weak equivalence. From the fact that the composite map $B\to \glob(\mathbf{S}^0) \to \vI$ is a trivial fibration and the \ttt, we deduce that the unique map of flows $\glob(\mathbf{S}^0) \to \vI$ is a weak equivalence. Consider the commutative diagram of flows
\[
\xymatrix@C=3em@R=3em
{
	\glob(\mathbf{S}^0) \ar@{=}[r] \fd{c_1} & \glob(\mathbf{S}^0) \fd{} \\
	\glob(\mathbf{D}^1) \fr{} & \cocartesien \vI
}
\]
The bottom horizontal map $\glob(\mathbf{D}^1)\to \vI$ is a weak equivalence, being a trivial q-fibration. By the \ttt, we deduce that $c_1:\glob(\mathbf{S}^0) \to \glob(\mathbf{D}^1)$ is a weak equivalence as well, and therefore a trivial cofibration since it is a q-cofibration which is one-to-one on states. The induction hypothesis is therefore proved for $n=1$. Suppose that the induction hypothesis is proved for $n\geq 1$. Using Proposition~\ref{connected-colim-glob} and the pushout diagram of spaces
\[
\xymatrix@C=3em@R=3em
{
	\mathbf{S}^{n-1} \fr{}\fd{} & \mathbf{D}^{n}\fd{}\\\
	\mathbf{D}^{n} \fr{} & \cocartesien\mathbf{S}^{n}
}
\]
we obtain the commutative diagram of flows 
\[
\xymatrix@C=3em@R=3em
{
	\glob(\mathbf{S}^{n-1}) \fr{c_n}\fd{c_n} & \glob(\mathbf{D}^{n})\fd{}\ar@{->}[rr] && \vI\ar@{=}[d]\\
	\glob(\mathbf{D}^{n}) \fr{} & \cocartesien\glob(\mathbf{S}^{n}) \fr{c_{n+1}}& \glob(\mathbf{D}^{n+1}) \fr{} & \vI
}
\]
Using the induction hypothesis, we deduce that the map $\glob(\mathbf{D}^{n})\to \glob(\mathbf{S}^{n})$ is a trivial cofibration, being a pushout of the trivial cofibration $c_n$. All maps $\glob(\mathbf{D}^{N})\to \vI$ for $N\geq 0$ are trivial q-fibrations, and hence trivial fibrations. Using the \ttt, we obtain the induction hypothesis for $n+1$. We have proved that all maps of $\cell(I^{gl}_{\geq 1})$ are trivial cofibrations. 

Now we can conclude the proof as follows. Let $X$ be a flow containing at least one execution path and let $X^{cof}$ be a q-cofibrant replacement of $X$. Consider the flow $\Mon(X^{cof})$ defined by the pushout diagram of flows
\[
\xymatrix@C=3em@R=3em
{
	X^0 \fr{} \fd{} & X^{cof} \fd{} \\
	\{0\} \fr{} & \cocartesien\Mon(X^{cof}).
}
\]
By Proposition~\ref{pseudo-leftproper}, the canonical map $X^{cof}\to \Mon(X^{cof})$ is a weak equivalence. Consequently, we can suppose without loss of generality that $X^0=\{0\}$ and that $X$ is a cellular object of the q-model structure of flows. Write the canonical map $\varnothing \to X$ as a composite $\varnothing \longrightarrow X^0 \longrightarrow X^1 \longrightarrow  X$ such that the map $X^0\to X^1$ belongs to $\cell(\{c_0\})$ and such that $X^1\to X$ belongs to $\cell(I^{gl}_{\geq 1})$. In particular, the map $X^1\to X$ is a trivial cofibration by the first part of the proof. Factor the canonical map $X^1\to \mathbf{1}$ as a composite $X^1\to X^\infty \to \mathbf{1}$ such that the left-hand map belongs to $\cell(I^{gl}_{\geq 1})$ and such that the right-hand map belongs to $\inj(I^{gl}_{\geq 1})$. It means that $X$ is weakly equivalent to $X^\infty$. Since the map $X^\infty \to \mathbf{1}$ is bijective on states, it is injective with respect to $C:\varnothing\to \{0\}$. Since, moreover, $X^\infty$ contains an execution path, it is also injective with respect to $c_0:\glob(\mathbf{S}^{-1})\subset \glob(\mathbf{D}^{0})$. Thus, the map $X^\infty \to \mathbf{1}$ is a weak equivalence, being a trivial fibration. We deduce that every flow in $(\C,\F,\W)$ is weakly equivalent to $\varnothing$, $\{0\}$ or $\mathbf{1}$. Since the full subcategory of $\dtop$ generated by the three objects $\varnothing$, $\{0\}$ and $\mathbf{1}$ is $\uu,$ the homotopy category of $(\C,\F,\W)$ is then a categorical localization of $\uu$. We deduce the inclusion $\overline{\W}_{DK}\subset \W$. The set of generating cofibrations $I^{gl}\cup \{C\}$ is tractable. Therefore, by \cite[Theorem~1.4]{henry2020minimal}, there exists a minimal model structure $(\C,\F,\W)$ with respect to the set of generating cofibrations $I^{gl}\cup \{C\}$. In this case, there is also the inclusion $\W\subset \overline{\W}_{DK}$ and the proof is complete since a model structure is characterized by its classes of cofibrations and weak equivalences.
\epf

\begin{cor} \label{DK-explicit1}
	The minimal model structure on flows with respect to the left-lifting of the cofibrations of the Ilias model structure has three homotopy types. 
\end{cor}

\bpf It is a consequence of Proposition~\ref{cof-left-induite} and Theorem~\ref{DK-explicit}.
\epf

\section{Right-lifting the q-model structure of flows}

We want to prove that the q-model structure of flows can be transferred along the right adjoint $\cat\subset \dtop$. At first, we recall:

\bth \label{GHKRS} (Kan-Quillen, see \cite[proof of Theorem~1 of Section~II.4 ]{MR36:6480} and \cite[Theorem~11.3.2]{ref_model2} or for an abstract presentation \cite[Theorem~2.2.1]{HKRS17}) Let $\mathcal{M}$ and $\mathcal{N}$ be two locally presentable categories. Let $(\C,\F,\W)$ be a combinatorial model structure of $\mathcal{M}$ such that all objects are fibrant. Consider a categorical adjunction $L:\mathcal{M} \dashv \mathcal{N}:U$. Suppose that there exists a factorization of the diagonal of $\mathcal{N}$ as a composite $X\stackrel{\tau}\to \cocyl(X) \stackrel{\pi} \to X\p X$ such that $U(\tau)$ is a weak equivalence of $\mathcal{M}$ and such that $U(\pi)$ is a fibration of $\mathcal{M}$ for all objects $X$ of $\mathcal{N}$. Then there exists a unique combinatorial model structure on $\mathcal{N}$ such that the class of fibrations is $U^{-1}(\F)$ and such that the class of weak equivalences is $U^{-1}(\W)$. If the set of generating (trivial resp.) cofibrations of $(\C,\F,\W)$ is $I$ ($J$ resp.), then the set of generating (trivial resp.) cofibrations of the model structure of $\mathcal{N}$ is $L(I)$ ($L(J)$ resp.).
\eth


In the terminology of this note, Theorem~\ref{three2} means that the q-model structure of flows has a right-lifting to the category of small topologically enriched categories which is minimal.

\bth \label{three2} There exists a unique model structure $(\cattop)_q=(\C_q,\F_q,\W_q)$ on $\cattop$ such that: 
\begin{itemize}[leftmargin=*]
	\item The set of generating cofibrations is $\{\II(\glob(\Sp^{n-1})) \subset \II(\glob(\Di^{n})) \mid n\geq 0\} \cup \{\II(C),\II(R)\}$.
	\item The set of generating trivial cofibrations is $\{\II(\glob(\Di^{n} \p \{0\})) \subset \II(\glob(\Di^{n}\p [0,1])) \mid n\geq 0\}$.
	\item A map of small enriched categories $f:X\to Y$ is a weak equivalence if and only if $\Obj(f):\Obj(X)\to \Obj(Y)$ is a bijection and for all $(\alpha,\beta)\in \Obj(X)\p \Obj(X)$, the continuous map $X(\alpha,\beta)\to X(f(\alpha),f(\beta))$ is a weak homotopy equivalence.
	\item A map of small enriched categories $f:X\to Y$ is a fibration if and only if for all $(\alpha,\beta)\in \Obj(X)\p \Obj(X)$, the continuous map $X(\alpha,\beta)\to X(f(\alpha),f(\beta))$ is a q-fibration of spaces.
\end{itemize}
Moreover, this model structure is minimal and all objects are fibrant. The left Quillen adjoint $\II:(\dtop)_q \to (\cat)_q$ is not a left Quillen equivalence.
\eth 

\bpf 
Consider the right adjoint $\cattop\subset \dtop$. Let $X$ be a small enriched category. Let $\cocyl(X)$ be the small enriched category having the same objects as $X$ and such that the space of morphisms $\cocyl(X)(\alpha,\beta)$ is equal to the topological space $\ttop([0,1],X(\alpha,\beta))$ with the continuous composition law defined for any triple $(\alpha,\beta,\gamma)$ of objects of $X$ as the composite: 
\begin{multline*}
\ttop([0,1],X(\alpha,\beta)) \p \ttop([0,1],X(\beta,\gamma)) \iso \ttop([0,1],X(\alpha,\beta) \p X(\beta,\gamma)) \\\longrightarrow \ttop([0,1],X(\alpha,\gamma)).
\end{multline*}
The composition law is clearly associative. The identity of $\cocyl(X)(\alpha,\alpha)$ (the space of morphisms in $\cocyl(X)$ from $\alpha$ to itself) is the constant map $\id_\alpha:[0,1] \to X(\alpha,\alpha)$. For all small enriched categories $X$, for all $(\alpha,\beta)\in \Obj(X)\p \Obj(X)$, the map $X(\alpha,\beta)\iso \ttop(\{0\},X(\alpha,\beta)) \to \ttop([0,1],X(\alpha,\beta))=\cocyl(X)(\alpha,\beta)$ is a trivial q-fibration of spaces and the map $\cocyl(X)(\alpha,\beta)=\ttop([0,1],X(\alpha,\beta)) \to \ttop(\{0,1\},X(\alpha,\beta)) \iso X(\alpha,\beta) \p X(\alpha,\beta)$ is a q-fibration of spaces. Using Theorem~\ref{GHKRS}, the q-model structure of $\dtop$ right induces a combinatorial model structure on $\cattop$. The model structure is minimal because it is combinatorial and all its objects are fibrant.

Let $X$ be an enriched small category. In $\dtop$, the map $X^{cof}\to X$ is a trivial q-fibration of flows. It means that for all $\alpha\in \Obj(X)$, $\P_{\alpha,\alpha}X^{cof}\to \P_{\alpha,\alpha}X$ is a trivial q-fibration of spaces. Therefore the map $\II(X^{cof})(\alpha,\alpha) = \{\id_{\alpha}\} \sqcup \P_{\alpha,\alpha}X^{cof} \to X(\alpha,\alpha)=\P_{\alpha,\alpha}X$ cannot be a weak homotopy equivalence. It implies that the map $\II(X^{cof}) \to X$ cannot be a weak equivalence of $(\cat)_q$. We deduce that the left Quillen adjoint $(\dtop)_q\to (\cat)_q$ is not homotopically surjective, and therefore that it is not a left Quillen equivalence.
\epf

We have $\II(\{0\}) \to (\{0\iso 1\})^{cof} \in (\C_{DK}\cap \W_{DK}) \backslash (\C_{q}\cap \W_{q})$. Thus, $\id:(\cattop)_{DK} \to (\cattop)_{q}$ cannot be a left Quillen adjoint. We have $R:\{0,1\}\to \{0\} \in \C_{q} \backslash \C_{DK}$. It implies that $\id:(\cattop)_{q} \to (\cattop)_{DK}$ cannot be a left Quillen adjoint either.

%

\end{document}